\documentclass[twoside]{article}[11pt]
\usepackage{a4}
\usepackage{verbatim}
\usepackage{amsmath}
\usepackage{amssymb}
\usepackage{amsthm}
\usepackage{definitions}
\usepackage[dvips]{graphicx}

\def\Inv{\mathcal{INV}}

\def\cDeg{\mathcal{DEG}}
\def\tr{\operatorname{tr}}
\def\Xdeg{X_{\operatorname{deg}}}
\def\Xinv{X_{\operatorname{inv}}}

\author{Jaros\l{}aw Buczy\'nski}
\date{6 July 2007}

\title{Some quasihomogeneous Legendrian varieties\footnote{
The article is a part of the research project 
N20103331/2715 funded by Polish financial 
means for science in years 2006-2008.
The article will be included in the author's PhD thesis.
Author's e-mail: jabu@mimuw.edu.pl
}
}
\begin{document}
\maketitle

%\section*{Things to correct}

%\textbf{
%}

\begin{abstract}
We construct a family of examples of Legendrian subvarieties in
 some projective spaces.
Although most of them are singular,
a new example of smooth Legendrian variety in dimension 8 is in this family.
The 8-fold has interesting properties:
it is a compactification of the special linear group,
a Fano manifold of index 5 and Picard number 1.
%A careful analysis of its properties leads us to conclusions     
%announced here and described precisely in some separate notes.

\end{abstract}

\tableofcontents

\section*{Acknowledgements}

The author is especially grateful to Sung Ho Wang for initiating the research. 
Also very special thanks to Insong Choe for his support and invitation
to KIAS (Korea Institute for Advanced Study)
where the author could (among many other attractions) meet Sung Ho Wang. 
The author acknowledges the help of Jaros\l{}aw Wi\'s{}niewski, Grzegorz Kapustka, 
Micha\l{} Kapustka, Michel Brion, Micha\l{} Krych, Joseph Landsberg, Laurent Manivel,
and Andrzej Weber.
Also thanks to an anonymous referee for pointing out some interesting
references, as well as for his other comments.

\section{Introduction}

Real Legendrian subvarieties are classical objects of differential
geometry and they have been investigated for ages.
However, complex Legendrian subvarieties in a projective space (see
\S\ref{notation_Legendrian} for the definition) are much
more rigid and only few smooth and compact examples were known
(see \cite{bryant}, \cite{landsbergmanivel04}, \cite{jabu_toric}):
\begin{enumerate}
  \item 
    linear subspaces;
  \item
    some homogeneous spaces called subadjoint varieties: 
    the product of a line and  a quadric $\P^1 \times Q^n$
    and five exceptional cases:
    \begin{itemize}
      \item
        twisted cubic curve  $\P^1 \subset \P^3$,
      \item
        Grassmannian $Gr_L(3,6) \subset \P^{13}$ of Lagrangian subspaces in $\C^6$,
      \item
        full Grassmannian  $Gr(3,6)\subset \P^{19}$,
      \item
        spinor variety $\mathbb{S}_6\subset \P^{31}$ 
         (i.e. the homogeneous $\mathbf{SO}(12)$-space 
         parametrising the vector subspaces of          dimension 6 contained 
         in a non-degenerate quadratic cone in $\C^{12}$) and 
      \item
        the 27-dimensional $E_7$-variety in $\P^{55}$ corresponding to the marked root: 
        \includegraphics[width=0.2\textwidth]{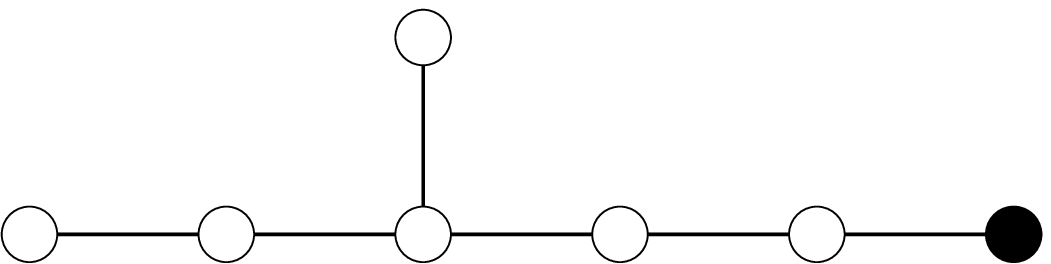};
    \end{itemize}
  \item
    every smooth projective curve admits a Legendrian embedding in $\P^3$  \cite{bryant};
  \item
    a family of surfaces birational to the Kummer $K3$-surfaces \cite{landsbergmanivel04};
  \item
    the blow up of $\P^2$ in three general points \cite{jabu_toric}. 
\end{enumerate}
In this article we present a new example in dimension 8 (see theorem
\ref{theorem_classify_invertible}(b)). Also we show how does the
construction generalise to give new examples in dimensions 5 and 14 
(see section \ref{section_other_examples}) and
finally we announce a result which will produce plenty of such examples
(see section \ref{section_hyperplane}).

 The original motivation for studying Legendrian subvarieties in a complex
projective space comes from the studies of contact Fano manifolds%
\footnote{
A complex projective manifold $M$ of dimension $2n+1$ is called \textbf{a contact manifold},
if there exists a vector subbundle $F \subset TM$ of rank $2n$,
such that the map $F \otimes F \lra TM/F$ determined by the Lie bracket is nowhere degenerate.
In such a case $F$ is called \textbf{a contact distribution}.
A projective manifold is \textbf{Fano}, if the anticanonical bundle is ample.
}
(see \cite{wisniewski}, \cite{kebekus_lines1}, \cite{4authors}):
The variety of tangent directions to the minimal rational curves through a fixed point on a contact Fano 
manifold make a Legendrian subvariety in the projectivisation of the fibre of contact distribution.
The adjoint varieties 
(i.e.~the closed orbit of the adjoint action of a simple Lie group $G$ on $\P(\gotg)$)
are the only known examples of contact Fano manifolds 
and they give rise to the homogeneous Legendrian varieties%
\footnote{
The groups of types $B$ and $D$ give rise to $\P^1 \times Q^n$. 
The five exceptional groups  $G_2$, $F_4$, $E_6$,  $E_7$, $E_8$ 
make the exceptional homogeneous Legendrian varieties. The groups of types $A$ and $C$ 
are somewhat special --- see \cite{landsbergmanivel04}, \cite{jabu06}.
}. 

From our considerations here, some other potential applications come into
the view --- see sections \ref{section_other_examples} and
\ref{section_compactification}.

\medskip

Before we present our results precisely in section \ref{results}, we must
introduce some notation. We need the notation of 
\S\ref{notation_begins}-\S\ref{notation_middle_end} to state the
results and also \S\ref{notation_middle_begin}-\S\ref{notation_ends} to prove them.

\subsection{Notation and definitions}\label{section_basics}

For this article we fix an integer $m\ge 2$.

\subsubsection{Vector space $V$}\label{notation_V}
\label{notation_begins}
Let $V$ be a vector space over complex numbers $\C$ of dimension $2m^2$,
 which we interpret as a space of pairs of $m \times m$ matrices. 
The coordinates are: $a_{ij}$ and $b_{ij}$ for $i,j \in \{1,\ldots m\}$. 
By $A$ we denote the matrix $(a_{ij})$ and similarly for $B$ and $(b_{ij})$.

\medskip

By $\P(V)$ we mean the naive projectivisation of $V$, i.e.~the
quotient  $(V \backslash \{0\})/ \C^*$.

\medskip

Given two $m\times m$ matrices $A$ and $B$, by $(A,B)$ we denote the point of the vector space $V$, 
while by $[A,B]$ we denote the point of the projective space $\P(V)$.

\medskip

Sometimes, we will represent some linear maps $V\lra V$ and
some 2-linear forms $V\otimes V \lra \C$ as $2m^2 \times 2m^2$ matrices. In such a case
we will assume the coordinates on $V$ come in the lexicographical
order:
$$
a_{11}, \ldots, a_{1m}, a_{21}, \ldots, a_{mm}, b_{11}, \ldots, b_{1m},
b_{21}, \ldots, b_{mm}.
$$

\subsubsection{Symplectic form $\omega$}
\label{notation_omega}
On $V$ we consider the standard symplectic form 
%$\omega:= \sum \ud a_{ij} \wedge \ud b_{ij}$  so that:
\begin{equation}
 \label{properties_of_inner_product}
\omega\big((A,B), (A', B')\big):= 
\sum_{i,j} (a_{ij}b'_{ij} -a'_{ij}b_{ij}) = \tr \left( A(B')^T - A'B^T \right).
\end{equation}
Further we set $J$ to be the
matrix of $\omega$:
$$
J:=
\left[ 
\begin{array}{cc}
0          & \Id_{m^2}\\
-\Id_{m^2} & 0
\end{array}
\right].
$$

\subsubsection{Lagrangian and Legendrian subvarieties}
\label{notation_Legendrian}

A linear subspace $W \subset V$ is called \emph{Lagrangian}
%(resp.~\emph{coisotropic})
 if the  $\omega$-perpen\-dicular subspace
$W^{\perp_\omega}$ is equal to %(resp.~contained in) 
$W$. 
Equivalently, $W$ is Lagrangian if and only if $\omega|_W \equiv 0$ and
$\dim W$ is maximal possible, i.e.~equal to $\half \dim V$.

\smallskip

A subvariety $Z \subset V$ is called \emph{Lagrangian} 
%(resp.~\emph{coisotropic})
if for every smooth point $z \in Z$ 
the tangent space $T_z Z \subset V$ is Lagrangian%
%(resp.~coisotropic)
.
In particular, if $Z$ is Lagrangian, then $\dim Z= \half \dim V= m^2$.

\smallskip

A subvariety $X \subset \P(V)$ is defined to be \emph{Legendrian} if its affine cone
$\hat{X}\subset V$ is Lagrangian.
In particular, if $X$ is Legendrian, then $\dim X= \half \dim V - 1= m^2-1$.

\subsubsection{Varieties $Y$, $\Xinv(m)$ and $\Xdeg(m,k)$}
\label{notation_Y}

We consider the following subvariety of $\P(V)$:
\begin{equation}
\label{defining_equations_of_Y}
Y:=\left\{[A,B]\in \P(V)  \mid AB^T = B^T A =
\lambda^2 \Id_m \textrm{ for some } \lambda \in \C\right\}.
\end{equation} 
The square at $\lambda$ seems to be irrelevant here, but it slightly simplifies the notation in the 
proofs  of theorem \ref{theorem_classify_invertible}(b)
and proposition \ref{proposition_orbits_of_G}(ii)

\label{notation_Xinv_Xdeg}

Further we define two types of subvarieties of $Y$:
$$
\Xinv(m):=\overline{
\bigg\{\left[g,\left(g^{-1}\right)^T\right] \in \P(V) \mid \det g = 1 \bigg\}}
$$
$$
\Xdeg(m,k):=\Big\{[A,B] \in \P(V) \mid  AB^T=B^TA=0, \ \rk A
  \le k, \ \rk B \le m-k  \Big\} 
$$
where $k \in {0,1,\ldots m}$.
The varieties $\Xdeg(m,k)$ have been also studied by \cite{strickland}
and \cite{mehta_trivedi}. $\Xinv(m)$ (especially $\Xinv(3)$) 
is the main object of this article.

\subsubsection{Automorphisms $\psi_{\mu}$} \label{notation_psi}

For any $\mu \in \C^*$ we let $\psi_{\mu}$ be the following linear
automorphism of $V$:
$$
\psi_{\mu}\big((A,B)\big) := (\mu A, \mu^{-1} B).
$$

Also the induced automorphism of $\P(V)$ will be denoted in the same
way:
$$
\psi_{\mu}\big([A,B]\big) := [\mu A, \mu^{-1} B].
$$

\medskip

The notation introduced so far is sufficient to state the results of this
paper (see section \ref{results}), but to prove them we need a few more
notions. 
%For the convenience of reader we group all the notation in
%one place.

\label{notation_middle_end}

\subsubsection{Groups $G$ and $\widetilde{G}$, Lie algebra $\gotg$ and
   their representation}
\label{notation_G}

\label{notation_middle_begin}

We set $\widetilde{G}:=\Gl_m \times \Gl_m$ and let it act on $V$ by:
$$
(g,h) \in \widetilde{G}, \ g,h \in \Gl_m, \ (A,B) \in V
$$
$$
(g,h) \cdot (A,B) := (g^T A h, g^{-1} B (h^{-1})^T).
$$
This action preserves the symplectic form $\omega$.

We will mostly consider the restricted action of 
$G:=\Sl_m \times \Sl_m < \widetilde{G}$.

We also set $\gotg:= \gotsl_m \times \gotsl_m$ to be the Lie algebra of $G$ and we have the
tangent action of $\gotg$ on $V$:
$$
(g,h) \cdot (A,B) = (g^T A + A h, - g B - B h^T).
$$
Though we denote the action of the groups $G$, $\widetilde{G}$ 
and the Lie algebra $\gotg$ by the same $\cdot$ 
we hope it will not lead to any confusion. 
Also the induced action of $G$ and $\widetilde{G}$ on $\P(V)$ will
be denoted by $\cdot$.

\subsubsection{Orbits $\Inv^m$ and $\cDeg^m_{k,l}$}
\label{notation_Inv_Deg}

We define the following sets:
$$
\Inv^m:=
\bigg\{\left[g,\left(g^{-1}\right)^T\right] \in \P(V) \mid \det g = 1 \bigg\},
$$
$$
\cDeg^m_{k,l}:= \Big\{[A,B] \in \P(V) \mid  AB^T=B^TA=0, \ \rk A = k, \
\rk B = l \Big\},
$$
so that $\Xinv(m)=\overline{\Inv^m}$ and $\Xdeg(m,k)= \overline{\cDeg^m_{k,m-k}}$. 

Clearly, if $k+l>m$ then $\cDeg^m_{k,l}$ is empty, so whenever
speaking of $\cDeg^m_{k,l}$ we will assume $k+l \le m$.

\subsubsection{Elementary matrices $E_{ij}$ and points $p_1$ and $p_2$}
\label{notation_Eij}

Let $E_{ij}$ be the elementary $m \times m$ matrix with
unit in the $i^{\textrm{th}}$ row and the $j^{\textrm{th}}$ column and zeroes
elsewhere.

\label{notation_p1_p2}

We distinguish two points 
$p_1 \in \cDeg^m_{1,0}$ and
$p_2 \in \cDeg^m_{0,1}$:
$$
p_1:=
\left[ E_{mm}, 0 \right]
\
\textrm{ and }
\
p_2:=
\left[ 0, E_{mm} \right]
$$

These points will be usually chosen as nice representatives of the
closed orbits
$\cDeg^m_{1,0}$ and $\cDeg^m_{0,1}$.

\subsubsection{Tangent cone}\label{properties_of_tangent_cone}

We recall the notion of the tangent cone and a few among many of its
properties. For more details and the proofs we refer to 
\cite[lecture 20]{harris} and \cite[III.\S3,\S4]{mumford}.

For an irreducible Noetherian scheme $X$ over $\C$ and a closed point $x\in X$ 
we consider the local ring $\ccO_{X,x}$ and we 
let $\gotm_x$ to be the maximal ideal in $\ccO_{X,x}$.
Let 
$$
R:=\bigoplus_{i=0}^{\infty} \left(\gotm_x^i \slash \gotm_x^{i+1} \right)
$$
where $\gotm_x^0$ is just the whole $\ccO_{X,x}$. 
Now we define \emph{the tangent cone $TC_x X$ at $x$ to $X$} to be $\Spec R$.

If $X$ is a subscheme of an affine space $\A^n$ (which we will usually
assume to be an affine piece of a projective space)
the tangent cone at $x$ to $X$ can be understood as a subscheme of $\A^n$. 
Its equations can be derived from the ideal of $X$.
For simplicity assume $x=0 \in \A^n$ and then the polynomials defining $TC_0 X$ 
are the lowest degree homogeneous parts of the polynomials in the ideal of $X$.

Another interesting point-wise definition is that $v \in TC_0 X$ is a closed point if and only if 
there exists a holomorphic map $\varphi_v$ from
the disc $D_t:=\{ t\in \C  \ : \ \lvert t \rvert\ < \delta \}$ to $X$,
such that $\varphi_v(0) = 0$
 and the first non-zero coefficient in the Taylor expansion in $t$ of $\varphi_v(t) $
is $v$, i.e.:
$$ \begin{array}{rccl}
\varphi_v: & D_t & \lra & X\\
           &  t  & \mapsto & t^k v + t^{k+1}v_{k+1} +\ldots
\end {array}$$

\smallskip

We list some of the properties of the tangent cone, that will be used
freely in the proofs:

\begin{itemize}
\item[(1)]
The dimension of every component of $TC_x X$ is equal to the dimension of $X$. 
\item[(2)] $TC_x X$ is naturally embedded in the Zariski tangent space to $X$ at $x$ and 
$TC_x X$ spans the tangent space.
\item[(3)] 
$X$ is regular at $x$ if and only if $TC_x X$ is equal (as scheme) to
	  the tangent space.
%\item[(4)]
%If $TC_x X$ is reduced, then $X$ is reduced at $x$.
\end{itemize}

%\begin{prf}
%The property (4) follows  from the definition and Krull Intersection Theorem 
%(see for example \cite[cor. 5.4]{eisenbud}).
%\end{prf}

\subsubsection{Submatrices - extracting rows and columns}\label{notation_submatrices}

Assume $A$ is an $m\times m$ matrix and $I,J$ are two sets of indices of
cardinality $k$ and $l$ respectively:
$$
I:=\{i_1,i_2, \ldots, i_k | 1 \le i_1 <  i_2 < \ldots < i_k \le m\},
$$
$$ 
J:=\{j_1,j_2, \ldots, j_l | 1 \le j_1 < j_2 < \ldots < j_l \le m\}.
$$
Then we denote by $A_{I,J}$ the $(m-k) \times (m-l)$ submatrix of $A$ obtained by 
removing rows of indices $I$ and columns of indices $J$.
Also for a set of indices $I$ we denote by $I'$ the set of $m-k$ indices complementary to $I$.

\smallskip

We will also use a simplified version of the above notation, when we
remove only a single column and single row:
$A_{ij}$ denotes the $(m-1) \times (m-1)$ submatrix of $A$ obtained by removing $i$-th row
and $j$-th column, i.e. $A_{ij} = A_{\{i\},\{j\}}$

\smallskip

Also in the simplest situation where we remove only the last row and
the last column, we simply write $A_m$, so that $A_m = A_{mm} = A_{\{m\},\{m\}}$.

\label{notation_ends}

\subsection{Main Results} \label{results}

In this note we give a classification%
\footnote{This problem was suggested by Sung Ho Wang.} 
of Legendrian subvarieties
in $\P(V)$ that are contained in $Y$.

\begin{theo}
\label{classification_theorem}
Let projective space $\P(V)$, varieties $Y$, $\Xinv(m)$,
 $\Xdeg(m,k)$ and automorphisms $\psi_{\mu}$ be defined as in
 \S\ref{notation_begins}-\S\ref{notation_middle_end}.
Assume $X\subset\P(V)$ is an irreducible subvariety.
Then $X$ is Legendrian and contained in $Y$ if and only if 
$X$ is one of the following varieties:
\begin{itemize}
\item[1.]
$X = \psi_{\mu}(\Xinv(m))$ for some $\mu \in \C^*$ or 
\item[2.]
$X=\Xdeg(m,k)$ for some $k \in \{0,1,\ldots m\}$.

\end{itemize}
\end{theo}

The idea of the proof of theorem \ref{classification_theorem} is based on
the observation that every Legendrian subvariety that is contained in $Y$
must be invariant under the action of group $G$. This is
explained in section \ref{section_action}. A proof of the theorem is
presented in section \ref{section_classification}.

Also we analyse which of the above varieties appearing in 1.~and 2.~are smooth:

\begin{theo}
    \label{theorem_classify_invertible}
  With the definition of $\Xinv(m)$ as in \S\ref{notation_Xinv_Xdeg},
  the family $\Xinv(m)$ contains the following varieties:
  \begin{itemize}
    \item[(a)] 
      $\Xinv(2)$ is a linear subspace.
    \item[(b)] 
      $\Xinv(3)$ is smooth, its  Picard group is generated by a hyperplane section.
      Moreover $\Xinv(3)$ is a compactification of $\Sl_3$ and it is
      isomorphic to a hyperplane section of Grassmannian $Gr(3,6)$. 
      The connected component of $Aut(\Xinv(3))$ is equal to $G=\Sl_3 \times \Sl_3$ 
      and $\Xinv(3)$ is not a homogeneous space.
    \item[(c)]
      $\Xinv(4)$ is the 15 dimensional spinor variety $\mathbb{S}_6$.
    \item[(d)] 
      For $m\ge 5$, the variety $\Xinv(m)$ is singular.
  \end{itemize}
\end{theo}

A proof of the theorem is explained in section \ref{section_invertible}.

Variety $\Xinv(3)$ is not yet described as a Legendrian subvariety
variety, so it is our new smooth example of dimension 8.

\begin{theo}
\label{theorem_smooth_degenerate}
With the definition of $\Xdeg(m)$ as in \S\ref{notation_Xinv_Xdeg},
variety $\Xdeg(m,k)$ is smooth if and only if $k=0$ , $k=m$ or $(m,k)=(2,1)$. 
In the first two cases, $\Xdeg(m,0)$ and $\Xdeg(m,m)$ are linear spaces, 
while $\Xdeg(2,1) \simeq \P^1 \times \P^1 \times \P^1 \subset \P^7$.
\end{theo}

A proof of the theorem is presented in section \ref{section_degenerate}.

\medskip

The results of theorems 
\ref{classification_theorem},
\ref{theorem_classify_invertible}
can be generalised in (at least) three different directions:

\subsection{Generalisation 1: Representation theory}
    \label{section_other_examples}

The interpretation of theorem
\ref{theorem_classify_invertible} (b) and (c) 
can be following:
We take the exceptional Legendrian variety $Gr(3,6)$, slice it with a linear section 
and we get a description, that generalised to matrices of bigger size gives 
the bigger exceptional Legendrian variety $\bS_6$.
Similar connection can be established between other exceptional Legendrian varieties.

\medskip

For instance, assume that $V^{sym}$ is a vector space of dimension $2 \binom{m+1}{2}$,
which we interpret as the space of pairs of $m \times m$ symmetric matrices $A,B$.
Now in $\P(V^{sym})$ consider the subvariety $\Xinv^{sym}(m)$, which is the closure of the following set:
$$
\{[A,A^{-1}]\in \P(V^{sym}) | A=A^T \textrm{ and } \det A=1\}.
$$
\begin{theo}
All the varieties $\Xinv^{sym}(m)$ are Legendrian and we have:
\begin{itemize}
\item[(a)] $\Xinv^{sym}(2)$ is a linear subspace.
\item[(b)] $\Xinv^{sym}(3)$ is smooth and it is isomorphic to a
	   hyperplane section of Lagrangian Grassmannian $Gr_L(3,6)$.
\item[(c)] $\Xinv^{sym}(4)$ is smooth and it is the 9 dimensional Grassmannian variety $Gr(3,6)$.
\item[(d)] For $m\ge 5$, the variety $\Xinv^{sym}(m)$ is singular.
\end{itemize}
\end{theo}

The proof goes exactly as the proof of theorem \ref{theorem_classify_invertible}.

\medskip

Similarly, we can take $V^{skew}$ to be a vector space of dimension $2 \binom{2m}{2}$,
which we interpret as the space of pairs of $2m \times 2m$ skew-symmetric matrices $A,B$.
Now in $\P(V^{skew})$ consider subvariety $\Xinv^{skew}(m)$, which is the closure of the following set:
$$
\{[A, - A^{-1}]\in \P(V^{skew}) | A= - A^T \textrm{ and } \pf A=1\}.
$$
\begin{theo}
All the varieties $\Xinv^{skew}(m)$ are Legendrian and we have:
\begin{itemize}
\item[(a)] $\Xinv^{skew}(2)$ is a linear subspace.
\item[(b)] $\Xinv^{skew}(3)$ is smooth and it is isomorphic to a hyperplane section of the spinor variety $\bS_6$.
\item[(c)] $\Xinv^{skew}(4)$ is smooth and it is the 27 dimensional $E_7$ variety.
\item[(d)] For $m\ge 5$, the variety $\Xinv^{skew}(m)$ is singular.
\end{itemize}
\end{theo}

Here the only difference is that we replace the determinants by the
Pfaffians of the appropriate submatrices
and also for the previous cases we will be picking some diagonal matrices
as nice representatives.
Since there is no non-zero skew-symmetric diagonal matrix,
we must modify a little bit our calculations, but there is no
essential difference in the technique.

\medskip

Neither $\Xinv^{sym}(3)$ nor $\Xinv^{skew}(3)$ 
have been described as smooth Legendrian subvarieties.

\medskip

Therefore we have established a new connection 
between the subadjoint varieties of the 4 exceptional
groups $F_4$, $E_6$, $E_7$ and $E_8$.
A similar connection was obtained by \cite{landsbergmanivel02}.

%Also it can be interesting to understand why this construction fails for the last step.
%For $Gr_L (3,6)$, $Gr(3,6)$, $\bS_6$, we just took a rank 4 Jordan algebra (corresponding respectively 
%to reals, complex numbers and quaternions) 
%and related to it the appropriate variety.
%So what happens if we take a $4\times 4$ octonionic hermitian matrices?
%One obstruction for defining the variety is that the determinants of degree 3 are not well defined
%(because octonions are not associative),
%so it is hard to say what $A^{-1}$ should be. 
%But note that to define $\Xinv^{sym}(4)$, $\Xinv(4)$ and  $\Xinv^{skew}(4)$, 
%we only use quadratic equations which can be expressed in a form of some matrix multiplication
%or as some combination of some $2 \times 2$ minors (or $4\times 4$ Pfaffians).
%So the only obstruction now is the noncommutativity 
%and somehow  it is not an obstruction for the quaternion case.

\subsection{Generalisation 2: Hyperplane section}\label{section_hyperplane}

The variety $\Xinv(3)$ is the first described example of smooth non-homogeneous Legendrian variety of dimension
bigger than $2$ (see \cite{bryant}, \cite{landsbergmanivel04}, \cite{jabu_toric}).
But this example is very close to a homogeneous one, namely is isomorphic to a hyperplane section of $Gr(3,6)$, 
which is a well known Legendrian variety. 
So a natural question arises, 
whether a general hyperplane section of other Legendrian varieties admits Legendrian embedding. 
The answer is yes and we explain it (as well as many conclusions from this surprisingly simple observation)
in \cite{jabu_hyperplane}.

\subsection{Generalisation 3: Group compactification}\label{section_compactification}

Theorem \ref{theorem_classify_invertible}(b) says that $\Xinv(3)$ is a
smooth compactification of $\Sl_3$. 
In \cite{jabu_compactifications} we study a generalisation of this construction 
(which is not really related to Legendrian varieties) 
to find a family of compactifications of $\Sl_n$, 
which contains the smooth compactification of $\Sl_3$ 
and can be easily smoothened (by a single blow up of a closed orbit) for
$n=4$.

\section{$G$-action and its orbits}\label{section_action}

In \cite{jabu06} we prove:

\begin{theo} \label{theorem_ideal_and_group}
  Let $X \subset \P(V)$ be a Legendrian subvariety (see \S\ref{notation_Legendrian} for definition).
  Consider the following map:
  $$
    H^0(\ccO_{\P(V)}(2)) \simeq \Sym^2 V^* \ni  q = (x \mapsto x^T M(q) x) 
    \stackrel{\rho}{\mapsto} 2J \cdot M(q) \in \gotsp(V).
  $$
  where $M(q)$ is the $(2m^2) \times (2m^2)$ matrix of $q$ and $J$ is the matrix of the symplectic   form as in \S\ref{notation_V}.
  Let $\I_2(X) \subset \Sym^2 V^*$ be the vector space of quadrics containing $X$. Then:
  \begin{itemize} 
    \item
      $\rho(\I_2(X))$ is a Lie subalgebra of $\gotsp(V)$ tangent to a closed subgroup 
      $$
        \overline{\exp\Big(\rho\big(\I_2(X)\big)\Big)} < \Sp(V).
      $$
    \item
      We have the natural action of $\Sp(V)$ on $\P(V)$. 
      The group 
      $\overline{\exp\Big(\rho\big(\I_2(X)\big)\Big)}$
      is the maximal connected subgroup in $\Sp(V)$ 
      which under this action preserves $X \subset \P(V)$.
  \end{itemize}
\end{theo}

\begin{prf}
  See \cite[cor.~4.4, cor.~5.5, lem.~5.6]{jabu06}.
\end{prf}

Recall the definition of $Y$ in \S\ref{notation_Y}.

The following polynomials are in 
the homogeneous ideal of $Y$ (below $i,j$ are indices that run through $\{1,\ldots, m\}$, $k$ is a summation index):
\begin{subequations}\label{equations_of_Y}
\begin{align}
\sum_{k=1}^m a_{ik} b_{ik}  - \sum_{k=1}^m a_{1k} b_{1k}
\label{equations_of_Y1}\\
\sum_{k=1}^m a_{ik} b_{jk} \textrm{ for $i \ne j$}
\label{equations_of_Y2}\\
\sum_{k=1}^m a_{ki} b_{ki}  - \sum_{k=1}^m a_{k1} b_{k1}
\label{equations_of_Y3}\\
\sum_{k=1}^m a_{ki} b_{kj} \textrm{ for $i \ne j$}
\label{equations_of_Y4}
\end{align}
\end{subequations}

These equations simply come from eliminating $\lambda$ from the
defining equation of $Y$ --- see equation \eqref{defining_equations_of_Y}.

For the statement and proof of the following proposition, recall our
notation of 
\S\ref{notation_V}, 
\S\ref{notation_omega}, 
\S\ref{notation_Legendrian}, 
\S\ref{notation_G} and
\S\ref{notation_Eij}.

\begin{prop}
\label{proposition_action_of_G}
Let $X\subset \P(V)$ be a Legendrian subvariety.
If $X$ is contained in $Y$ 
then $X$ is preserved by the induced action of $G$ on $\P(V)$.
\end{prop}

\begin{prf}
Let $\I_2(X)$ be as in the theorem \ref{theorem_ideal_and_group} 
and define  $\I_2(Y)$ analogously.
Clearly $\I_2(Y) \subset \I_2(X)$. 
By theorem \ref{theorem_ideal_and_group} it is enough to calculate that 
$\gotg \subset \rho\left(\I_2(Y)\right)$ or that the images of the quadrics 
\eqref{equations_of_Y1}--\eqref{equations_of_Y4} under $\rho$ generate $\gotg$.

We deal in details of the proof only for $m=2$. 
There is no difference between this case and the general
one, except for the complexity of notation%
%(which should be 3-dimensional at this point to be clear enough).
.

Let us take the quadric
$$
q_{ij}:=\sum_{k=1}^m a_{ik} b_{jk} = a_{i1} b_{j1} + a_{i2} b_{j2} 
$$
for any $i,j \in \{1,\ldots, m\} = \{1,2\}$.
Also let $Q_{ij}$ be the $2m^2 \times 2m^2$ symmetric matrix
corresponding to $q_{ij}$. For instance:
$$
Q_{12} =
\left[ 
\begin{array}{cccccccc}
0&0 &0&0 &0&0 &\half & 0 \\
0&0 &0&0 &0&0 &0 & \half \\
0&0 &0&0 &0&0 &0 & 0 \\
0&0 &0&0 &0&0 &0 & 0 \\
0&0 &0&0 &0&0 &0 & 0 \\
0&0 &0&0 &0&0 &0 & 0 \\
\half&0 &0&0 &0&0 &0 & 0 \\
0&\half &0&0 &0&0 &0 & 0
\end{array}
\right] .
$$

So choose an arbitrary $(A,B) \in V$ and at the
moment we want to think of it as of a single vertical $2 m^2$-vector:
$(A,B) = [a_{11}, a_{12},a_{21}, a_{22},b_{11}, b_{12},b_{21},
b_{22}]^T$, so that the following multiplication makes sense:
$$
\rho(q_{12}) = 
2 J \cdot Q_{12} \cdot (A,B) = 
$$
$$
= \left[ 
\begin{array}{cccccccc}
0&0 &0&0 &1&0 &0 & 0 \\
0&0 &0&0 &0&1 &0 & 0 \\
0&0 &0&0 &0&0 &1 & 0 \\
0&0 &0&0 &0&0 &0 & 1 \\
-1&0 &0&0 &0&0 &0 & 0 \\
0&-1 &0&0 &0&0 &0 & 0 \\
0&0 &-1&0 &0&0 &0 & 0 \\
0&0 &0&-1 &0&0 &0 & 0
\end{array}
\right]
\left[ 
\begin{array}{cccccccc}
0&0 &0&0 &0&0 &1 & 0 \\
0&0 &0&0 &0&0 &0 & 1 \\
0&0 &0&0 &0&0 &0 & 0 \\
0&0 &0&0 &0&0 &0 & 0 \\
0&0 &0&0 &0&0 &0 & 0 \\
0&0 &0&0 &0&0 &0 & 0 \\
1&0 &0&0 &0&0 &0 & 0 \\
0&1 &0&0 &0&0 &0 & 0
\end{array}
\right]
\left[ 
\begin{array}{c}
a_{11}\\
a_{12}\\
a_{21}\\
a_{22}\\
b_{11}\\
b_{12}\\
b_{21}\\
b_{22}
\end{array}
\right]
=
$$
$$
=
\left[ 
\begin{array}{cccccccc}
0&0 &0&0 &0&0 &0 & 0 \\
0&0 &0&0 &0&0 &0 & 0 \\
1&0 &0&0 &0&0 &0 & 0 \\
0&1 &0&0 &0&0 &0 & 0 \\
0&0 &0&0 &0&0 &-1 & 0 \\
0&0 &0&0 &0&0 &0 & -1 \\
0&0 &0&0 &0&0 &0 & 0 \\
0&0 &0&0 &0&0 &0 & 0 
\end{array}
\right]
\left[ 
\begin{array}{c}
a_{11}\\
a_{12}\\
a_{21}\\
a_{22}\\
b_{11}\\
b_{12}\\
b_{21}\\
b_{22}
\end{array}
\right]
=
$$
$$
=
\left[ 
\begin{array}{c}
0\\
0\\
a_{11}\\
a_{12}\\
-b_{21}\\
-b_{22}\\
0\\
0
\end{array}
\right]
\stackrel{\textrm{back to the matrix notation} }{=}
\left(
\left[ 
\begin{array}{cc}
0& 0 \\
a_{11} & a_{12}
\end{array}
\right] , \
\left[ 
\begin{array}{cc}
-b_{21} & -b_{22}\\
0& 0
\end{array}
\right]
\right)
= 
$$
$$
=
\left(
\left[ 
\begin{array}{cc}
0 &1 \\
0 &0
\end{array}
\right]^T
\left[ 
\begin{array}{cc}
a_{11} & a_{12}\\
a_{21} & a_{22}\\
\end{array}
\right] , \ 
- \left[ 
\begin{array}{cc}
0& 1 \\
0 &0
\end{array}
\right]
\left[ 
\begin{array}{cc}
b_{11} & b_{12}\\
b_{21} & b_{22}
\end{array}
\right]
\right)
 \ = \ (E_{12}^T A, \  -E_{12} B)
$$

Going along exactly the same calculations, we see that:

%\begin{equation} \label{sl_elementary_action}
$$
2 J \cdot Q_{ij} \cdot (A,B)  \ = \ (E_{ij}^T A, \  -E_{ij} B)
$$
%\end{equation}

Next in the ideal of $Y$ we have the following quadrics: $q_{ij}$ for
$i \ne j$ (see \eqref{equations_of_Y2}) and $q_{ii}-q_{11}$ (see
\eqref{equations_of_Y1}).
By taking images under $\rho$ of the linear combinations of those quadrics we can get an arbitrary
traceless matrix $g\in \gotsl_m$ acting on $V$ in the following way:
$$
g \cdot (A,B) = (g^T A, -g B).
$$
Exponentiate this action of $\gotsl_m$ to get the action of $\Sl_m$:
$$
g \cdot (A,B) = (g^T A, g^{-1} B).
$$

This proves of that the action of subgroup $\Sl_m \times 0 < G = \Sl_m \times \Sl_m$
indeed preserves $X$ as claimed in the lemma. 
The action of the other component $0 \times \Sl_m$ is calculated in the same way, but
using quadrics \eqref{equations_of_Y3}--\eqref{equations_of_Y4}.
\end{prf}

\subsection{Invariant subsets}\label{section_orbits}

Recall our notation of 
\S\ref{notation_V},
\S\ref{notation_Y},
\S\ref{notation_psi},
\S\ref{notation_G} and
\S\ref{notation_Inv_Deg}.
Here we want to decompose $Y$ into a union of some $G$-invariant
subsets, most of which are orbits.

\begin{prop}
\label{proposition_orbits_of_G}
\hfill
\begin{itemize}
\item[(i)]
The sets $\Inv^m$, $\psi_{\mu}(\Inv^m)$ and $\cDeg^m_{k,l}$ 
are $G$-invariant and they are all contained in $Y$.
\item[(ii)]
$Y$ is equal to the union of all $\psi_{\mu}(\Inv^m)$ (for
  $\mu \in \C^*$) and all $\cDeg^m_{k,l}$ (for integers $k,l\ge
  0$, $k+l\le m$).
\item[(iii)]
Every $\psi_{\mu}(\Inv^m)$ is an orbit of the action of $G$. 
If $m$ is odd, then $\Inv^m$ is isomorphic (as algebraic variety) to
$\Sl_m$. Otherwise if $m$ is even, then $\Inv^m$ is isomorphic to 
$(\Sl_m / \Z_2) $. In both cases, $\dim \psi_{\mu}(\Inv^m) =  \dim
\Inv^m = m^2 -1$. 
\end{itemize}
\end{prop}

\begin{prf}
The proof of part (i) is an explicit verification from the definitions
in \S\ref{section_basics}.

\medskip

To prove part (ii), assume $[A,B]$ is a point of $Y$, so $A B^T=B^TA=\lambda^2\Id_m$. 
First assume that the ranks of both matrices are maximal: 
$$
\rk A = \rk B =m. 
$$
Then $\lambda$ must be non-zero and $B=\lambda^2 (A^{-1})^T$.
Let $d:=(\det A)^{-\frac{1}{m}} $ so that 
$$
\det (dA)=1
$$
and let 
$\mu:= \frac{1}{d \lambda} $.
Then we have:
$$
[A,B] =
 \left[A,\lambda^2 \left(A^{-1}\right)^T\right] =
 \left[ \frac{dA}{d \lambda }, d\lambda \left((dA)^{-1}\right)^T\right] =
$$
$$
% = [ dA, \mu^{-2} ((dA)^{-1})^T]
= \left[ \mu (dA), \mu^{-1} \left((dA)^{-1}\right)^T\right]= 
 \psi_{\mu} \left(\left[ (dA), \left((dA)^{-1}\right)^T\right]\right).
$$
Therefore $[A,B] \in \psi_{\mu}(\Inv^m)$.

Next, if either of the ranks is not maximal:
$$
\rk A < m \textrm{ or } \rk B < m
$$
then by \eqref{defining_equations_of_Y} we must have $AB^T = B^TA =
0$. 
So $[A,B] \in \cDeg^m_{k,l}$ for $k=\rk A$ and $l=\rk B$. 

\medskip

Now we prove (iii).
The action of $G$ commutes with $\psi_{\mu}$:
$$
(g,h)\cdot \psi_{\mu}\big([A,B]\big) = \psi_{\mu}\big((g,h) \cdot [A,B]\big).
$$
So to prove $\psi_{\mu}(\Inv^m)$ is an orbit it is enough to prove
that $\Inv^m$ is an orbit, which follows from the definitions of the
action and $\Inv^m$.

We have the following epimorphic map:
$$
\begin{array}{rcl}
\Sl_m &\lra & \Inv^m\\
g &\mapsto& [g,(g^{-1})^T]
\end{array}
$$
If $[g_1,(g_1^{-1})^T] = [g_2,(g_2^{-1})^T]$ then we must have
$g_1 = \alpha g_2$ and $g_1 = \alpha^{-1} g_2$ for some 
$\alpha \in \C^*$. 
Hence $\alpha^2=1$ and $g_1 = \pm g_2$. If $m$ is odd and 
$g_1\in \Sl_m$ then $-g_1 \notin \Sl_m$ so $g_1 = g_2$. So $\Inv^m$ is 
either isomorphic to $\Sl_m$ or to $\Sl_m /\Z_2$ as stated.
\end{prf}

From proposition \ref{proposition_orbits_of_G}(ii) we conclude that 
$\Xinv(m)$ is an equivariant compactification of $\Sl_m$ (if $m$ is odd) 
or $\Sl_m / \Z_2$ (if $m$ is even).
See \cite{timashev} and references therein 
for the theory of equivariant compactifications. 
In the setup of \cite[\S8]{timashev}, 
this is the compactification corresponding to the representation $W \oplus W^*$,
where $W$ is the standard representation of $\Sl_m$. 
Therefore some properties of $\Xinv(m)$ could also be read from the general description 
of group compactifications.

\begin{prop}
\label{degenerate_orbits_of_G}
\hfill
\begin{itemize}
\item[(i)]
The dimension of $\cDeg^m_{k,l}$ is $(k+l)(2m-k-l)-1$.
In particular, if $k+l=m$ then the dimension is equal to $m^2-1$.
\item[(ii)]
$\cDeg^m_{k,l}$ is an orbit of the action of $G$, unless $m$ is even
  and $k=l= \half m$. 
\item[(iii)]
If $m\ge 3$, then there are exactly two closed orbits of the action of
$G$: $\cDeg^m_{1,0}$ and $\cDeg^m_{0,1}$.
\end{itemize}
\end{prop}

\begin{prf}
Part (i) follows from \cite[prop 2.10]{strickland}.

\medskip

For part (ii) 
let $[A,B] \in \cDeg^m_{k,l}$ be any point.
By Gauss elimination and elementary linear algebra,
we can prove that there exists 
$(g,h) \in G$ such that $[A',B']:=(g,h) \cdot [A,B]$ is a pair of
diagonal matrices. Moreover, if $k+l < m$ then we can choose $g$ and
$h$ such that:
$$
A':=\diag(\underbrace{1,\ldots, 1}_{k},\underbrace{0,\ldots, 0}_{l},
\underbrace{0,\ldots, 0}_{m-k-l}),
$$\nopagebreak
$$
B':=\diag(\underbrace{0,\ldots, 0}_{k},\underbrace{1,\ldots, 1}_{l},
\underbrace{0,\ldots, 0}_{m-k-l}).
$$
Hence $\cDeg^m_{k,l} = G\cdot [A',B']$  and this finishes the proof in the case
$k+l < m$. 

So assume $k+l=m$. Then we can choose $(g,h)$ such that:
$$
A':=\diag(\underbrace{1,\ldots, 1}_{k},\underbrace{0,\ldots, 0}_{l}),
$$
$$
B':=\diag(\underbrace{0,\ldots, 0}_{k},\underbrace{d,\ldots, d}_{l}),
$$
for some $d \in \C^*$. 
If $k\ne l$, then set $e:=d^{\frac{1}{l-k}}$ and let
$$
g':=\diag(\underbrace{e^l,\ldots, e^l}_{k},\underbrace{e^{-k},\ldots, e^{-k}}_{l}).
$$
Clearly $\det (g')=1$ and:
$$
(g', \Id_m) \cdot [A',B'] = 
\left[
\diag(\underbrace{e^l,\ldots, e^l}_{k},\underbrace{0,\ldots,0}_{l}),
\diag(\underbrace{0,\ldots, 0}_{k},\underbrace{d e^k,\ldots,d e^k}_{l})
\right]
$$
where
$$
d e^k = d^{1+\frac{k}{l-k}} = d^{\frac{l}{l-k}} = e^l.
$$
So rescaling we get:
$$
(g', \Id_m) \cdot [A',B'] = 
\left[
\diag(\underbrace{1,\ldots, 1}_{k},\underbrace{0,\ldots,0}_{l}),
\diag(\underbrace{0,\ldots, 0}_{k},\underbrace{1,\ldots,1}_{l})
\right]
$$
and this finishes the proof of (ii).

\medskip

For part (iii),
denote by $W_1$ (respectively, $W_2$) the standard representation of 
the first (respectively, the second) component of $G= \Sl_m \times \Sl_m$.
Then our representation $V$ is isomorphic to $(W_1 \otimes W_2) \oplus (W_1^*\otimes W_2^*)$. 
For $m\ge 3$ the representation $W_i$ is not isomorphic to $W_i^*$ and
therefore $V$ is a union of two irreducible non-isomorphic
representations, so there are exactly two closed orbits of this
action on $\P(V)$.  These orbits are simply $\cDeg^m_{1,0}$
and $\cDeg^m_{0,1}$.
\end{prf}

\subsection{Action of $\widetilde{G}$}

Recall the notation of 
\S\ref{notation_V},
\S\ref{notation_G} and
\S\ref{notation_Inv_Deg}.

The action of $\widetilde{G}$ extends the action of $G$, but it does
not preserve $\Xinv(m)$. 
So we will only consider the action of $\widetilde{G}$ when speaking of
$\Xdeg(m,k)$.

We have properties analogous to proposition
\ref{degenerate_orbits_of_G} (ii) and  (iii) but with no exceptional cases:

\begin{prop}
\label{properties_of_Gtilde}
\hfill
\begin{itemize}
\item[(i)]
Every $\cDeg^m_{k,l}$ is an orbit of the action of $\widetilde{G}$.
\item[(ii)]
For every  $m$ there are exactly two closed orbits of the action of
$\widetilde{G}$: $\cDeg^m_{1,0}$ and $\cDeg^m_{0,1}$.
\end{itemize}
\end{prop}

\begin{prf}
It goes exactly as the proof of proposition
\ref{degenerate_orbits_of_G} (ii) and (iii).
\end{prf}

\section{Legendrian varieties in $Y$}\label{main_part}

In this section we prove the main results of the article.

\subsection{Classification}\label{section_classification}

We start with proving the theorem \ref{classification_theorem}. 
For this we use our notation of section \ref{section_basics}.

\begin{prf}
First assume $X$ is Legendrian and contained in $Y$.
If $X$ contains a point $[A,B]$ where both
 $A$ and $B$ are invertible,  then by proposition
\ref{proposition_action_of_G}
 it must contain the orbit of $[A,B]$, which by proposition 
\ref{proposition_orbits_of_G}(ii) and (iii) is equal to 
$\psi_{\mu}(\Inv^m)$ for some $\mu \in \C^*$. But dimension of $X$ is $m^2-1$ which is
 exactly the dimension of $\psi_{\mu}(\Inv^m)$ (see proposition
 \ref{proposition_orbits_of_G}(iii)), 
so
$$
X =\overline{\psi_{\mu}(\Inv^m)} = \psi_{\mu}(\Xinv(m)).
$$

\medskip

On the other hand if $X$ does not contain any point $[A,B]$ where both $A$ and $B$ are invertible then 
in fact $X$ is contained in the locus $Y_0:=\{[A,B]:AB^T=B^T A=0\}$. 
This locus is just the union of all $\cDeg^m_{k,l}$ and its irreducible components are
the closures of $\cDeg^m_{k,m-k}$, which are exactly $\Xdeg(m,k)$.
 So in particular every irreducible component has
 dimension $m^2-1$ (see proposition \ref{degenerate_orbits_of_G}(i)) 
and hence $X$ must be one of these components.

\medskip

Therefore it remains to show that all these varieties are Legendrian. 

The fact that $\Xdeg(m,k)$ is a Legendrian variety follows from
\cite[pp524--525]{strickland}. Strickland proves there that the affine
cone over $\Xdeg(m,k)$ (or $W(k, m-k)$ in the notation of
\cite{strickland}) is the closure of a conormal bundle. Conormal bundles
are classical examples of Lagrangian varieties.

\medskip

Since $\psi_{\mu}$ preserves the symplectic form $\omega$, it is
enough to prove that $\Xinv(m)$ is Legendrian.

The group $G$ acts symplectically on $V$ and the action has an open
orbit on $\Xinv(m)$ --- see proposition \ref{proposition_orbits_of_G} (iii).
Thus the tangent spaces to the affine cone over $\Xinv(m)$ are Lagrangian if
 and only if just one tangent space at a point of the open orbit is
 Lagrangian.

So we take $[A,B]:=[\Id_m, \Id_m]$.
Now the affine tangent space to $\Xinv(m)$ at $[\Id_m,\Id_m]$ 
is the linear subspace of $V$ spanned by $(\Id_m,\Id_m)$ 
and the image of the tangent action of the Lie algebra $\gotg$. 
We must prove that for every four traceless matrices
$g,h, g', h'$ we have:
\begin{equation}
\label{equation_for_Inv_Legendrian1}
\omega \left((g,h) \cdot (\Id_m,\Id_m), \ (g',h') \cdot
(\Id_m,\Id_m)\right) = 0 \textrm{ and}
\end{equation}
\begin{equation}
\label{equation_for_Inv_Legendrian2}
\omega \left((\Id_m,\Id_m),  \ (g,h)\cdot (\Id_m,\Id_m)\right) = 0
\end{equation}

Equality \eqref{equation_for_Inv_Legendrian1} is true without the
assumption that the matrices have trace $0$:
$$
\omega \big((g,h) \cdot (\Id_m,\Id_m),  \ (g',h') \cdot
(\Id_m,\Id_m)\big) \ = 
$$\nopagebreak
$$
= \ \omega\Big( \left( g^T+h, \ -(g+h^T)\right), \quad
\left((g')^T+h', \ -(g' +(h')^T)\right)\Big) \ =
$$\nopagebreak
$$
\stackrel{\textrm{by \eqref{properties_of_inner_product}}}{=}
 \   \tr \Big( -\left(g^T+h\right) \left((g')^T + h'\right)   \  +
 \  \left(g+h^T\right) \left( g'+(h')^T\right)   \Big)  \ =
$$\nopagebreak
$$
= \ 0.
$$

For equality \eqref{equation_for_Inv_Legendrian2} we calculate:
$$
\omega \big((\Id_m,\Id_m),  \ (g,h) \cdot (\Id_m,\Id_m)\big) \ = 
$$\nopagebreak
$$
= \ \omega\Big((\Id_m,\Id_m) , \quad
\left(g^T+h, \ -(g +h^T)\right)\Big) \ =
$$\nopagebreak
$$
\stackrel{\textrm{by \eqref{properties_of_inner_product}}}{=}
- \tr (g^T+h)  - \tr (g+h^T) \ = \ 0.
$$

Hence we have proved that the closure of $\Inv^m$ is Legendrian.
\end{prf}

\subsection{Degenerate matrices}\label{section_degenerate}

Recall our notation of 
\S\ref{notation_V}, 
\S\ref{notation_Xinv_Xdeg},
\S\ref{notation_G},
\S\ref{notation_p1_p2} and
\S\ref{properties_of_tangent_cone}.

By \cite[prop.~1.3]{strickland} the ideal of $\Xdeg(m,k)$ is generated by
the coefficients of $AB^T$, the coefficients of $B^TA$, 
the $(k+1)\times(k+1)$-minors of $A$ and the $(m-k+1)\times(m-k+1)$-minors of $B$.
In short we will say that the equations of $\Xdeg(m,k)$ are given by:

\begin{equation}\label{equations_of_Xdeg}
AB^T=0, \ B^TA=0, \quad \rk(A) \le k, \ \rk(B) \le m-k.
\end{equation}

\begin{lem}
Assume $m\ge 2$ and $1\le k \le m-1$. Then:
\begin{itemize}
\item[(i)]
The tangent cone to $\Xdeg(m,k)$ at $p_1$ is a product of a linear space of dimension 
$(2m-2)$ and the affine cone over $\Xdeg(m-1,k-1)$.
\item[(i')]
The tangent cone to $\Xdeg(m,k)$ at $p_2$ is a product of a linear space of dimension 
$(2m-2)$ and the affine cone $\Xdeg(m-1,k)$.
\item[(ii)]
$\Xdeg(m,k)$ is smooth at $p_1$ if and only if $k=1$.
\item[(ii')]
$\Xdeg(m,k)$ is smooth at $p_2$ if and only if $k=m-1$.
\end{itemize}
\end{lem}

\begin{prf}
We only prove (i) and (ii), while (i') and (ii') follow in the same way by exchanging $a_{ij}$ and $b_{ij}$.
Consider equations \eqref{equations_of_Xdeg} of $\Xdeg(m,k)$ 
restricted to the affine neighbourhood of $p_1$ obtained by substituting $a_{mm}=1$. 
Taking the lowest degree part of these equations we get some of the equations of the tangent cone at $p_1$
(recall our convention on the notation of submatrices  --- see 
\S\ref{notation_submatrices}):
$$
b_{im}=b_{mi}=0, \ A_m B_m^T =0, \ B_m^T A_m = 0, 
$$
$$ 
\rk A_m \le k-1, \rk B_m \le m-k.
$$
These equations define the product of the linear subspace  $A_m=B_m=0, b_{im}=b_{mi}=0$ 
and the affine cone over $\Xdeg(m-1,k-1)$ embedded in the set of those pairs of matrices,
 whose last row and column are zero: $a_{im}=a_{mi}=0, b_{im}=b_{mi}=0$. 
So the variety defined by those equations is irreducible and its dimension is equal to
$(m-1)^2 + 2m-2 = m^2-1 = \dim \Xdeg(m,k)$. 
Since it contains the tangent cone we are interested in
 and by \S\ref{properties_of_tangent_cone}(1),
they must coincide as claimed in (i).

\medskip

Next (ii) follows immediately, since for $k=1$ the equations above reduce to 
$$
b_{im}=b_{mi}=0, \ \textrm{ and } A_m=0
$$
and hence the tangent cone is just the tangent space, so  $p_1$ is a smooth point of $\Xdeg(m,1)$.
Conversely,  if $k>1$ then $\Xdeg(m-1, k-1)$ is not a linear space, so by (i) 
the tangent cone is not a linear space either and $X$ is singular at $p_1$ --- see 
\S\ref{properties_of_tangent_cone}(3).
\end{prf}

Now we can prove theorem \ref{theorem_smooth_degenerate}:

\begin{prf}
It is obvious from the definition of $\Xdeg(m,k)$, 
that $\Xdeg(m,0)= \{A=0\}$ and $\Xdeg(m,m) = \{B=0\}$,
so these are indeed linear spaces.

Therefore assume $1\le k \le m-1$. 
But $\Xdeg(m,k)$ is
$\widetilde{G}$ invariant (see proposition \ref{properties_of_Gtilde}(i)) and so is its singular locus $S$. 
Hence $\Xdeg(m,k)$ is singular if and only if $S$ contains a closed orbit of $\widetilde{G}$.

So $\Xdeg(m,k)$ is smooth, 
if and only if it is smooth at both $p_1$ and $p_2$
(see proposition \ref{properties_of_Gtilde}(ii)), which 
(by lemma (ii) and (ii')) holds if and only if $k=1$ and $m=2$.

\smallskip
 
To finish the proof, 
it remains to verify what kind of variety is $\Xdeg(2,1)$.
Consider the following map:
$$
\P^1 \times \P^1 \times \P^1 \lra \P(V) \simeq \P^7
$$
$$
[\mu_1, \mu_2],[\nu_1,\nu_2],[\xi_1,\xi_2] \longmapsto
\left[
\xi_1
\left(
\begin{array}{cc}
\mu_1\nu_1 & \mu_1 \nu_2\\
\mu_2\nu_1 & \mu_2 \nu_2
\end{array}
\right),
\xi_2
\left(
\begin{array}{cc}
\mu_2\nu_2 & -\mu_2 \nu_1\\
-\mu_1\nu_2 & \mu_1 \nu_1
\end{array}
\right)
\right]
$$
Clearly this is a Segre embedding in appropriate coordinates.
The image of this embedding is contained in $\Xdeg(2,1)$ 
(see equation \eqref{equations_of_Xdeg}) 
and since dimension of $\Xdeg(2,1)$ 
is equal to the dimension of  $\P^1 \times \P^1 \times \P^1$ 
we conclude the above map gives an isomorphism of 
$\Xdeg(2,1)$ and  $\P^1 \times \P^1 \times \P^1$.

%It can be easily and explicitely seen, 
%that it is $\P^1 \times \P^1 \times \P^1$, 
%but let us just observe the following properties:
%Some quadratic equations generate the ideal of $\Xdeg(2,1)$ 
%(see \eqref{equations_of_Xdeg}), 
%$\Xdeg(2,1)$ is Legendrian by theorem \ref{classification_theorem}
%and it is smooth irreducible of dimension 3 
%(see proposition \ref{degenerate_orbits_of_G}(i)),
%so by \cite[thm 5.11]{jabu06} it must be exactly the product of $\P^1$'s.
\end{prf}

\subsection{Invertible matrices}\label{section_invertible}

Recall the notation of 
\S\ref{notation_V},
\S\ref{notation_Xinv_Xdeg},
\S\ref{notation_G},
\S\ref{notation_Inv_Deg} and
\S\ref{properties_of_tangent_cone}.

We wish to determine some of the equations of $\Xinv(m)$.
Clearly the equations of $Y$ (see \eqref{equations_of_Y}) 
are quadratic equations of $\Xinv(m)$.
To find other equations, we recall, that 
$$
\Xinv(m):=\overline{
\bigg\{\left[g,\left(g^{-1}\right)^T\right] \in \P(V) \mid \det g = 1 \bigg\}}
$$
But for a matrix $g$ with determinant 1 we know that the entries of
$(g^{-1})^T$
consist of the appropriate minors (up to sign)
of $g$.
Therefore we get many inhomogeneous equations satisfied by every pair
$\left(g, (g^{-1})^T\right) \in V$
(recall our convention on the notation of submatrices  --- see 
\S\ref{notation_submatrices}):
$$
\det (A_{ij}) = (-1)^{i+j} b_{ij} \  \textrm{ and } \  a_{kl} = (-1)^{k+l} \det(B_{kl})  
$$

To make them homogeneous, multiply two such equations appropriately:
\begin{equation}\label{minors_1_of_X}
\det(A_{ij}) a_{kl} = (-1)^{i+j+k+l} b_{ij} \det(B_{kl}).
\end{equation}
These are degree $m$ equations, 
which are satisfied by the points of $\Xinv(m)$ 
and we state the following theorem:

\begin{theo}
\label{theorem_Xinv3_is_smooth}
Let $m=3$. Then the quadratic equations \eqref{equations_of_Y1}--\eqref{equations_of_Y4} and
the cubic equations \eqref{minors_1_of_X} generate the ideal of 
$\Xinv(3)$. Moreover $\Xinv(3)$ is smooth.
\end{theo}

\begin{prf}
It is enough to prove that the scheme $X$ defined by equations
\eqref{equations_of_Y1}--\eqref{equations_of_Y4} 
and  \eqref{minors_1_of_X} is smooth, because the  reduced subscheme
of $X$ coincides with $\Xinv(3)$.

The scheme $X$ is $G$ invariant, 
hence as in the proof of theorem \ref{theorem_smooth_degenerate} and by
proposition \ref{degenerate_orbits_of_G}(iii) it is
 enough to verify smoothness 
at $p_1$ and $p_2$.
Since we have the additional symmetry here (exchanging $a_{ij}$'s with $b_{ij}$'s)
it is enough to verify the smoothness at $p_1$.

Now we calculate the tangent space to $X$ at $p_1$ by taking 
linear parts of the equations evaluated at $a_{33}=1$. 
From \eqref{equations_of_Y} we get that 
$$
b_{31}=b_{32}=b_{33}=b_{23}=b_{13}=0.
$$
Now from equations \eqref{minors_1_of_X} for $k=l=3$ and $i,j \ne 3$  we get the following evaluated equations:
$$
a_{i'j'} - a_{i'3}a_{3j'} = \pm b_{ij} B_{33} 
$$
(where $i'$ is either $1$ or $2$, which ever is different than $i$ and analogously for $j'$)
so the linear part is just $a_{i'j'}=0$. Hence by varying $i$ and $j$ we can get 
$$
a_{11}= a_{12}=a_{21}=a_{22}=0.
$$
Therefore the tangent space has codimension at least $9$, which is exactly the codimension of 
$\Xinv(3)$ --- see \ref{proposition_orbits_of_G}(iii). Hence $X$ is
smooth (in particular reduced) and $X=\Xinv(3)$. 
\end{prf}

To describe $\Xinv(m)$ for $m > 3$ we must find more equations. 

There is a more general version of the above  property of an inverse
of a matrix with determinant 1, which is less popular.

\begin{prop}\label{proposition_minors_2}
\hfill
\begin{itemize}
\item[(i)]
Assume $A$ is a $m \times m$ matrix of determinant $1$ and $I,J$ are
two sets of indices, both of cardinality $k$
(again recall our convention on indices and submatrices  --- see
section \ref{notation_submatrices}). Denote by $B:=(A^{-1})^T$. 
Then the appropriate minors are equal (up to sign):
$$
\det A_{I,J} = (-1)^{\Sigma I + \Sigma J} \det B_{I', J'}.
$$
\item[(ii)]
The coordinate free way to express these equalities is following: 
if  $W$ is a vector space of dimension $m$ 
and $f$ is a linear automorphism of $W$,
let $\Wedge{k} f$ be the induced automorphism of $\Wedge{k} W$.
If $\Wedge{m} f = \Id_{\Wedge{m} W}$  then:
$$
\Wedge{m-k} f  = \Wedge{k} \left( \Wedge{m-1} f\right) .
$$
\item[(iii)]
Consider the induced action of $G$ on the polynomials on $V$.
Then the vector space spanned by the set of equations of (i) for a fixed
	  $k$ is $G$ invariant.
\end{itemize}
\end{prop}

\begin{prf}
Part (ii) follows explicitly from (i), since 
if $A$ is a matrix of $f$, then the terms of  the matrices of  the maps 
$\Wedge{m-k} f$ and $\Wedge{k} ( \Wedge{m-1} f)$
 are exactly the appropriate minors of $A$ and $B$.

Part (iii) follows easily from (ii).

As for (i), 
%compare with a slightly similar statement in \cite[\S II.4.3]{sikorski}.
%?? what is a slightly similar statement???
%
we only sketch the proof, leaving the details to the reader and his or her linear algebra students.
Firstly, reduce to the case when $I$ and $J$ are just $\{1, \ldots k\}$ and the determinant of $A$
is possibly $\pm 1$ (which is where the sign shows up in the equality).
Secondly if both determinants $\det A_{I,J}$ and $\det B_{I', J'}$ are 
zero, then the equality is clearly satisfied. Otherwise assume for example $\det A_{I,J} \ne 0$.
Then performing the appropriate row and column operations we can change $A_{I,J}$ into a diagonal matrix, 
$A_{I', J}$ and $A_{I, J'}$ into the zero matrices 
and all these operations can be done without changing $B_{I',J'}$ 
nor $\det A_{I,J}$. 
Then the statement follows easily.
\end{prf}

In particular we get:

\begin{cor}
\label{corollary_equations_Xinv}
Assume $k$, $I$ and $J$ are as in proposition \ref{proposition_minors_2}(i).
\begin{itemize}
\item[(a)]
If $m$ is even and $k=\half m$, then the equation
$$
\det A_{I,J} = (-1)^{\Sigma I + \Sigma J} \det B_{I', J'}
$$
is homogeneous of degree $\half m$ and it is satisfied by points of $\Xinv(m)$.
\item[(b)]
If $0 \le k < \half m $ and $l= \half m - k$, then
$$
\left(\det A_{I,J}\right)^2 = \left(\det B_{I',J'}\right)^2 \cdot (a_{11}b_{11} + \ldots  a_{1m}b_{1m})^l
$$
is a homogeneous equation of degree $2(m-k)$ satisfied by points of $\Xinv(m)$.
\end{itemize}
\end{cor}

\begin{prf}
Clearly both equations are homogeneous.
If $\det A =1$ and $B=(A^{-1})^T$ then the following equations are satisfied: 
\begin{equation}\label{equation_for_cor1}
\det A_{I,J} = (-1)^{\Sigma I + \Sigma J} \det B_{I',J'},
\end{equation}
\begin{equation}\label{equation_for_cor2}
1 = (a_{11}b_{11} + \ldots  a_{1m}b_{1m})^l 
\end{equation}
(equation \eqref{equation_for_cor1} follows from proposition \ref{proposition_minors_2}(i)
and \eqref{equation_for_cor2} follows from $AB^T=\Id_m$).
Equation in (b) is just \eqref{equation_for_cor1} squared multiplied side-wise by \eqref{equation_for_cor2}.

So both equations in (a) and (b) are satisfied by every pair 
$\left( A, (A^{-1})^T\right)$ 
and by homogeneity also by
$\left(\lambda A, \lambda (A^{-1})^T \right)$. 
Hence (a) and (b) hold on an open dense subset of $\Xinv(m)$, so also on whole $\Xinv(m)$.
\end{prf}

We know enough equations of $\Xinv(m)$ to prove the theorem \ref{theorem_classify_invertible}:

\subsubsection{Case $m=2$ --- linear subspace}

\begin{prf}
To prove (a) just take the linear equations from  proposition \ref{proposition_minors_2}(i) for $k=1$:
$$
a_{ij} = \pm b_{i'j'}
$$
where $\{i,i'\} = \{j,j'\} = \{1,2\}$.
\end{prf}

\subsubsection{Case $m=3$ --- hyperplane section of $Gr(3,6)$}

\begin{prf}
For (b), $\Xinv(3)$ is smooth by theorem \ref{theorem_Xinv3_is_smooth} 
and it is a compactification of $\Inv^3 \simeq \Sl_3$ by proposition
\ref{proposition_orbits_of_G}(i) and (iii).

\paragraph{Picard group of $\Xinv(3)$.}

The complement of the open orbit 
$$
D:=\Xinv(3) \backslash \Inv^3
$$
must be a union of some orbits of $G$, each of them must have
dimension smaller than $\dim \Inv^3 =8$. So by propositions
\ref{proposition_orbits_of_G}(ii), (iii), \ref{degenerate_orbits_of_G}
(i) and (ii) the only candidates are $\cDeg^3_{1,1}$, $\cDeg^3_{0,1}$
and $\cDeg^3_{1,0}$.
We claim they are all contained in $\Xinv(3)$. It is enough to prove
that $\cDeg^3_{1,1} \subset \Xinv(3)$, since the other orbits are in
the closure of $\cDeg^3_{1,1}$. Take the curve in $\Xinv(3)$
parametrised by:
$$
\left[
\left(
\begin{array}{ccc}
t  &0  &0\\
0  &1  &0\\
0  &0  &t^{-1}\\
\end{array}
\right),
\left(
\begin{array}{ccc}
t^{-1} &0  &0\\
0      &1  &0\\
0      &0  &t\\
\end{array}
\right)
\right].
$$
For $t=0$ the curve meets $\cDeg^3_{1,1}$, which finishes the proof of
the claim.

Since $\dim \cDeg^3_{1,1} = 7$  (see proposition
\ref{degenerate_orbits_of_G}(i)), $D$ is a prime divisor.
We have $\Pic(\Sl_3) = 0$ and by \cite[prop. II.6.5(c)]{hartshorne} the Picard group of 
 $\Xinv(3)$ is isomorphic to $\Z$ with the ample generator $[D]$.

Next we check that $D$ is linearly equivalent 
(as a divisor on $\Xinv(3)$) to a hyperplane section $H$ of $\Xinv(3)$. 
Since we already know that $\Pic(\Xinv(3))=\Z \cdot [D]$, we must have
$H \stackrel{lin}{\sim} k D$ for some positive integer $k$.
But there are lines contained in $\Xinv(3)$ (for example those
contained in $\cDeg^3_{1,0} \simeq \P^2 \times \P^2$)\footnote{
Actually, the reader could also easily find explicitly some lines (or
even planes) which intersect the open orbit 
and conclude that $\Xinv(3)$ is covered by lines.
}.
So let $L\subset \Xinv(3)$ be any line and we intersect:
$$
D\cdot L = \frac{1}{k} H \cdot L = \frac{1}{k}.
$$
But the result must be an integer, so $k=1$ as claimed.

\paragraph{Complete embedding.}

Since $D$ itself is definitely not a hyperplane section of $\Xinv(3)$, 
the conclusion is that the Legendrian embedding of $\Xinv(3)$ is not given by a complete linear system.
The natural guess for a better embedding is the following:
%
%$$
%\Xinv(m):=\overline{
%\bigg\{\left[g,\left(g^{-1}\right)^T\right] \in \P(V) \mid \det g = 1 \bigg\}}
%$$
%
$$
X' := \overline{
\bigg\{\left[1,g,\Wedge{2} g\right] \in \P^{18} = \P(\C \oplus V)
\mid \det g = 1 \bigg\}}.
$$
(we note that $\Wedge{2} g = (g^{-1})^T$ for $g$ with $\det g=1$)
and one can verify that the projection from the point $[1,0,0]\in \P^{18}$
restricted to $X'$ gives an isomorphism with $\Xinv(3)$.

The Grassmannian $Gr(3,6)$ in its Pl\"u{}cker embedding can be described as the closure of:
$$
\bigg\{\left[1,g,\Wedge{2} g, \Wedge{3} g\right] \in \P^{19} = \P(\C \oplus V \oplus \C)
\mid g \in M_{3 \times 3}\bigg\}
$$
and we immediately identify $X'$ as the section $H:=\left\{ \Wedge{3} g=1 \right\}$ of the Grassmannian. 

Though it is not essential, we note that $H^1(\ccO_{Gr(3,6)}) = 0$ 
(see Kodaira vanishing theorem \cite[thm 4.2.1]{lazarsfeld}) and hence the above embedding of $\Xinv(3)$ 
is given by the complete linear system.

\paragraph{Automorphism group.}

It remains to calculate 
$Aut\left(\Xinv(3)\right)^0$ 
--- the connected component of the automorphism group.

The tangent Lie algebra of the group of automorphisms of a complex projective 
manifold is equal to the global sections of the tangent bundle, see \cite{akhiezer}.
A vector field on $\Xinv(3)$ is also a section of $T Gr(3,6)|_{\Xinv(3)}$ 
and we have the following short exact sequence:
$$
0 \lra T Gr(3,6)(-1) \lra T Gr(3,6) \lra T Gr(3,6)|_{\Xinv(3)} \lra 0
$$
The homogeneous vector bundle $T Gr(3,6)(-1)$ is isomorphic to 
$U^* \otimes Q \otimes \Wedge{3} U$, where $U$ is the universal subbundle in  $Gr(3,6) \times \C^6$ 
and $Q$ is the universal quotient bundle.
This bundle corresponds to an irreducible module of the parabolic subgroup in $\Sl_6$.
Calculating explicitly its highest weight and applying Bott formula 
\cite{ottaviani} we get that $H^1\big(T Gr(3,6)(-1)\big) = 0$. 
Hence every section of $T\Xinv(3)$ extends to a section of $TGr(3,6)$.
In other words, if $P < Aut(Gr(3,6)) \simeq \P\Gl_6$ 
is the subgroup preserving $\Xinv(3) \subset Gr(3,6)$,
then the restriction map $P \lra Aut\left(\Xinv(3)\right)^0$ is epimorphic.

The action of $\Sl_6$ on $\Wedge{3} \C^6$ preserves 
the natural symplectic form $\omega'$:
$$
\omega' : \Wedge{2}\left(\Wedge{3} \C^6\right) \lra  \Wedge{6} \C^6 \simeq \C.
$$
Since the action of $P$ on $\P\left(\Wedge{3} \C^6\right)$ 
preserves the hyperplane $H$ containing $\Xinv(3)$, 
it must also preserve $H^{\perp_{\omega'}}$, 
i.e.~$P$ preserves $[1,0,0,1]\in \P^{19} = \P(\C\oplus V \oplus \C)$.
Therefore $P$ acts on the quotient $H/(H^{\perp_{\omega'}}) = V$ 
and hence the restriction map factorises:
$$
P \lra Aut(\P(V), \Xinv(3))^0 \epi Aut(\Xinv(3))^0.
$$

By \cite{jabu_toric}, group $Aut(\P(V), \Xinv(3))^0$ 
is contained in the image of $\Sp(V) \lra \P\Gl(V)$,
so by theorem \ref{theorem_ideal_and_group},
proposition \ref{proposition_action_of_G}
and theorem \ref{theorem_Xinv3_is_smooth}
\[
Aut\left(\P(V), \Xinv(3)\right)^0= G.
\]
In particular $\Xinv(3)$ cannot be homogeneous as it contains 
more than one orbit of the connected component of automorphism group.

\end{prf}

We note that the  fact that $\Xinv(3)$ is not homogeneous 
can be also proved without calculating the automorphism group.
Since $\Pic \Xinv(3) \simeq \Z$, 
it follows from \cite[thm.~11]{landsbergmanivel04},
that $\Xinv(3)$ could only be one of the subadjoint varieties.
But none of them has $\Pic \simeq \Z$ and dimension 8.

\subsubsection{Case m=4 --- spinor variety $\mathbb{S}_6$}

\begin{prf}
To prove (c) we only need to take 30 quadratic equations of $Y$ as in
\eqref{equations_of_Y} and 36 quadratic 
equations from  corollary \ref{corollary_equations_Xinv} (a).
By proposition \ref{proposition_minors_2}(iii)  
the scheme $X$ defined by those quadratic equations is $G$ invariant.  
As in the proofs of theorems  \ref{theorem_smooth_degenerate} and \ref{theorem_Xinv3_is_smooth},
we only check that $X$ is smooth at $p_1$ and $p_2$ and 
conclude it is smooth everywhere, hence those equations indeed define $\Xinv(4)$.
 
Therefore $\Xinv(4)$ is smooth, irreducible and its ideal is  generated by quadrics, 
so  it falls into the classification of \cite[thm. 5.11]{jabu06}. 
Hence we have two choices for $\Xinv(4)$, whose dimension is $15$:
the product of a line and a quadric $\P^1\times Q_{14}$ or the spinor variety $\mathbb{S}_6$.
The homogeneous ideal of polynomials vanishing on $\P^1\times Q_{14}\subset \P^{31}$ 
is generated by $\dim (\Sl_2\times \mathbf{SO}_{16})=123$ 
linearly independent quadratic polynomials 
(see theorem \ref{theorem_ideal_and_group},
alternatively, one can calculate the equations explicitly --- see \cite[\S7.2]{jabu_arxiv} ). 
So $\Xinv(4)$, which by the above argument is generated by only 66 quadratic equations,
must be isomorphic to $\mathbb{S}_6$.
\end{prf}

\subsubsection{Case m$\ge 5$ --- singular varieties} 

\begin{prf}
Finally we prove (d).
We want to prove, that for $m\ge 5$ variety $\Xinv(m)$ is singular at $p_1$. 
To do that, we calculate the reduced tangent cone 
$$
T:=\big(TC_{p_1} \Xinv(m)\big)_{red}.
$$ 
From equations \eqref{equations_of_Y} we easily get the following 
linear and quadratic equations of $T$
(again we suggest to have a look at \S\ref{notation_submatrices}):
$$
b_{im} = b_{mi}=0, \quad A_m B_m^T = B_m^T A_m = \lambda^2 \Id_{m-1}
$$
for every $i\in \{1,\ldots m\}$  and some $\lambda\in \C^*$.

Next assume $I$ and $J$  are two sets of indices both of cardinality $k = \left\lfloor \half m \right\rfloor$
and such that neither $I$ nor $J$ contains $m$. 
Consider the equation of $\Xinv(m)$ as in corollary \ref{corollary_equations_Xinv}(b):
$$
\left(\det A_{I,J}\right)^2 = \left(\det B_{I',J'}\right)^2 \cdot (a_{11}b_{11} + \ldots  a_{1m}b_{1m})^l.
$$
To get an equation of $T$, we evaluate at $a_{mm}=1$ and take the lowest degree part,
which is simply 
$\left(\det \left((A_m)_{I,J}\right)\right)^2=0$.
Since $T$ is reduced, by varying $I$ and $J$
we get that: 
$$
\rk A_m \le m - 1 - k - 1 = \left\lceil \half m \right\rceil - 2 
$$
and therefore also:
$$
 A_m B_m^T = B_m^T A_m = 0.
$$

Hence $T$ is contained in the product of the linear space 
$W:=\{A_m=0, B=0\}$ and the affine cone $\hat{U}$ over the union of $\Xdeg(m-1,k)$ for 
$k\le \left\lceil \half m \right\rceil - 2$. We claim that $T= W \times \hat{U}$.
By proposition \ref{degenerate_orbits_of_G}(i), every component of $W \times \hat{U}$ 
has dimension $2m-2+(m-1)^2 = m^2-1=\dim \Xinv(m)$,
so by \S\ref{properties_of_tangent_cone}(1) the tangent cone must be a union of some of the components.
Therefore to prove the claim it is enough to find for every $k \le \left\lceil \half m \right\rceil - 2$ 
a single element of $\cDeg^{m-1}_{k, m-k-1}$ that is contained in the tangent cone. 

So take $\alpha$ and $\beta$ to be two strictly  positive integers such that 
$$
\alpha =  \left( \half m - k - 1 \right) \beta
$$
and consider the curve in $\P(V)$ with the following parametrisation:
$$
\left[
\diag\{
\underbrace{t^{\alpha},\ldots, t^{\alpha}}_{k},
\underbrace{t^{\alpha+\beta},\ldots,t^{\alpha+\beta}}_{m-k-1}, 1
\},
\diag\{
\underbrace{t^{\alpha+\beta},\ldots, t^{\alpha+\beta}}_{k},
\underbrace{t^{\alpha},\ldots,t^{\alpha}}_{m-k-1}, t^{2\alpha + \beta}
\}
\right].
$$
It is easy to verify that this family is contained in $\Inv^m$ for $t\ne 0$
and as $t$ converges to $0$, it gives rise to a tangent vector 
(i.e.~an element of the reduced tangent cone - see point-wise definition in
\S\ref{properties_of_tangent_cone}) that belongs to $\cDeg^{m-1}_{k, m-k-1}$.

So indeed $T= W \times \hat{U}$, which for $m\ge 5$ contains more than 1 component, 
hence cannot be a linear space. 
Therefore by \S\ref{properties_of_tangent_cone}(3) variety $\Xinv(m)$ is singular at $p_1$.
\end{prf}

\begin{rem}
Note that in both  cases of $\Xdeg(m,k)$ and $\Xinv(m)$, the reduced tangent cone 
is a Lagrangian subvariety in the fibre of the contact distribution.
This is not accidental as explained in \cite{jabu_tangent_cone}.
\end{rem}

\bibliography{references}

\bibliographystyle{alpha}

\end{document}